\documentclass{amsart}
\usepackage{amsmath}
\usepackage{amssymb}
\usepackage[english]{babel}
\usepackage{graphics}
\usepackage[all]{xy}
\usepackage{amsrefs}
\usepackage{amsmath}
\usepackage{amsthm}
\usepackage{amssymb}
\usepackage{amsfonts}
\usepackage{amsxtra}
\usepackage{tikz}
\usepackage{multicol}
\usepackage{subfigure}
\usepackage{centernot}
\usepackage{soul}
\usepackage{xcolor}
\usepackage{longtable}
\usepackage{layout}
\usepackage{geometry}
\usepackage{amsthm}

\usepackage{amstext} % for \text macro
\usepackage{array}   % for \newcolumntype macro
\newcolumntype{L}{>{$}c<{$}} % math-mode version of "l" column type
\usepackage{longtable}

\usepackage[colorlinks=true,linkcolor=blue,urlcolor=blue, citecolor=blue]{hyperref}
\usepackage{fullpage}
\usepackage{epsfig}
\usepackage{verbatim}
\usepackage{enumitem}
\usepackage{algorithmicx}
\usepackage{float}
\usepackage{algpseudocode}

\usepackage[justification=centering]{caption}

\newtheorem{Theorem}{Theorem}[section]

\newtheorem{Lemma}[Theorem]{Lemma}

\theoremstyle{definition}

\theoremstyle{remark}

\newcommand{\cod}{\operatorname{cod}}
\newcommand{\irr}{\operatorname{Irr}}
\newcommand{\PSL}{\normalfont{\mbox{PSL}}}
\newcommand{\PSU}{\normalfont{\mbox{PSU}}}

\usepackage{graphicx}
\usepackage{adjustbox}

\title{On the Characterization of Sporadic Simple Groups by Codegrees}

\author[M.~Dolorfino]{Mallory Dolorfino}
\author[L. ~Martin]{Luke Martin}
\author[Z. ~Slonim]{Zachary Slonim}
\author[Y. ~Sun]{Yuxuan Sun}
\author[Y. ~Yang]{Yong Yang}
\address{Mallory~Dolorfino\\ Kalamazoo College \\ Kalamazoo, Michigan, USA \\
  \href{mailto:mallory.dolorfino19@kzoo.edu}
  {{\ttfamily\upshape mallory.dolorfino19@kzoo.edu}}}
\address{Luke~Martin\\ Gonzaga University\\ Spokane, Washington, USA \\
  \href{mailto:lwmartin2019@gmail.com}
  {{\ttfamily\upshape lwmartin2019@gmail.com}}}
\address{Zachary~Slonim\\ University of California, Berkeley \\ Berkeley, California, USA \\
  \href{mailto:zachslonim@berkeley.edu}
  {{\ttfamily\upshape zachslonim@berkeley.edu}}}
\address{Yuxuan~Sun\\ Haverford College \\ Haverford, Pennsylvania, USA \\
  \href{mailto:ysun1@haverford.edu}
  {{\ttfamily\upshape ysun1@haverford.edu}}}
\address{Yong~Yang\\ Texas State University \\ San Marcos, Texas, USA \\
  \href{mailto:yang@txstate.edu}
  {{\ttfamily\upshape yang@txstate.edu}}}


\begin{document}

\maketitle

\begin{abstract}

    Let $G$ be a finite group and $\mathrm{Irr}(G)$ the set of all irreducible complex characters of $G$.  Define the codegree of $\chi \in \mathrm{Irr}(G)$ as $\mathrm{cod}(\chi):=\frac{|G:\mathrm{ker}(\chi) |}{\chi(1)}$ and denote by $\mathrm{cod}(G):=\{\mathrm{cod}(\chi)|\chi\in \mathrm{Irr}(G)\}$ the codegree set of $G$. Let $H$ be one of the $26$ sporadic simple groups. In this paper, we show that $H$ is determined up to isomorphism by cod$(H)$.

\end{abstract}

\begin{section}{Introduction}

Let $G$ be a finite group and $\mathrm{Irr}(G)$ the set of all irreducible complex characters of $G$. For any $\chi \in \mathrm{Irr}(G),$ define the codegree of $\chi,$ denoted by $\mathrm{cod}(\chi),$ as $\mathrm{cod}(\chi) := \frac{|G:\mathrm{ker}(\chi)|}{\chi(1)}.$ Then define the codegree set of $G$ as $\mathrm{cod}(G) := \{\mathrm{cod}(\chi)  |  \chi \in \mathrm{Irr}(G)\}.$ The concept of codegrees was originally considered in \cite{Chillag}, where the codegree was defined as $\mathrm{cod}(\chi) := \frac{|G|}{\chi(1)},$ and it was later modified to its current definition by \cite{Qian} so that $\mathrm{cod}(\chi)$ is the same for $G$ and $G/N$ when $N \leq \mathrm{ker}(\chi).$ Several properties of codegrees have been studied, such as the relationship between the codegrees and element orders, codegrees of $p$-groups, and groups with few codegrees.

The codegree set of a group is closely related to the character degree set of a group, defined as $\mathrm{cd}(G) := \{\chi(1) |  \chi \in \mathrm{Irr}(G)\}.$ The relationship between the character degree set and the group structure is an active area of research, and many properties of group structure are largely determined by the character degree set. In 1990, Bertram Huppert studied the character degrees of simple groups and made the following conjecture about the relationship between a simple group $H$ and a finite group $G$ that have equal character degree sets.

{\bf Huppert's Conjecture:} Let $H$ be a finite nonabelian simple group and $G$ a finite group such that $\mathrm{cd}(H) = \mathrm{cd}(G).$ Then $G \cong H \times A,$ where $A$ is an abelian group.

Huppert's conjecture has since been verified for many cases such as the alternating groups, sporadic groups, and simple groups of Lie type with low rank, but it has yet to be verified for simple groups of Lie type with high rank. Recently, a similar conjecture related to codegrees has been posed.

{\bf Codegree Version of Huppert's Conjecture:} Let $H$ be a finite nonabelian simple group and $G$ a finite group such that $\mathrm{cod}(H) = \mathrm{cod}(G).$ Then $G \cong H.$

This conjecture appears in the \emph{Kourovka Notebook of Unsolved Problems in Group Theory} as question 20.79 \cite{Khukhro}. It has been verified for $\mathrm{PSL}(2,q)$, $\mathrm{PSL}(3,4),$ $\mathrm{Alt}_7,$ $\mathrm{J}_1$, $^2B_2(2^{2f+1})$ where $f \geq 1$, $\mathrm{M}_{11}, \mathrm{M}_{12}, \mathrm{M}_{22}, \mathrm{M}_{23}$ and $\mathrm{PSL}(3,3)$ by \cites{Ahanjideh,Bahri,Gintz}. The conjecture has also been verified for $\PSL(3,q)$ and $\PSU(3,q)$ in \cites{lmy} and ${}^2G_2(q)$ in \cites{gzy}. Most of these results concern simple groups with less than $21$ character degrees \cite{Aziziheris}. However, in this paper, we provide a general proof verifying this conjecture for all sporadic simple groups. The methods used may be generalized to simple groups of Lie type, giving promising results for characterizing all simple groups by their codegree sets.

\begin{Theorem}\label{thm11}

    Let $H$ be a sporadic simple group and $G$ a finite group. If $\mathrm{cod}(G)= \mathrm{cod}(H)$, then $G \cong H$.

\end{Theorem}

Throughout the paper, we follow the notation used in Isaac's Book \cite{Isaacs} and the ATLAS of Finite Groups \cite{Atlas}.

\end{section}

\begin{section}{Preliminary Results}\label{theory}

We first reproduce some lemmas which will be used in later proofs.

\begin{Lemma}\label{oneextend}
    \cite{Moreto}*{Lemma 4.2} Let $S$ be a finite nonabelian simple group. Then there exists $1_S\neq \chi \in \mathrm{Irr}(S)$ that extends to $\mathrm{Aut}(S)$.
\end{Lemma}
\begin{Lemma}\label{prodsimple}
    \cite{James}*{Theorem 4.3.34} Let $N$ be a minimal normal subgroup of $G$ such that $N=S_1\times\dots\times S_t$ where $S_i \cong S$ is a nonabelian simple group for each $i=1,\dots,t$. If $\chi \in \mathrm{Irr}(S)$ extends to $\mathrm{Aut}(S)$, then $\chi\times\dots\times\chi\in\mathrm{Irr}(N)$ extends to $G$.
\end{Lemma}
\begin{Lemma}\label{perfect}
    \cite{Gintz}*{Remark 2.6} Let $G$ be a finite group and $H$ a finite nonabelian simple group with $\mathrm{cod}(G)=\mathrm{cod}(H)$. Then $G$ is a perfect group.
\end{Lemma}

\begin{Lemma}\label{order}
    \cite{Hung} Let $G$ be a finite group and $S$ a finite nonabelian simple group such that $\mathrm{cod}(S) \subseteq$ $\mathrm{cod}(G)$. Then $|S|$ divides $|G|$.
\end{Lemma}

\begin{Lemma}\label{codegreeSubset}
    Let $G$ be a finite group with $N \trianglelefteq G.$ Then $\cod(G/N) \subseteq \cod(G)$.
\end{Lemma}
\begin{proof}
From \cite{Isaacs}*{Lemma 2.22}, we can define $\irr(G/N) = \{ \hat{\chi}(gN) = \chi(g) \mid \chi \in \irr(G) \text{ and } N \subseteq \operatorname{ker}(\chi)\}$. Take any $\hat{\chi} \in \mathrm{Irr}(G/N)$. By definition, we know that $\hat{\chi}(1) = \chi(1),$ so the denominators of $\mathrm{cod}(\hat{\chi})$ and $\mathrm{cod}(\chi)$ are equal. In addition, $\mathrm{ker}(\hat{\chi}) \cong \mathrm{ker}(\chi) / N$, so $|\mathrm{ker}(\chi)| = |N| \cdot |\mathrm{ker}(\hat{\chi})|$. Thus $|G/N : \ker(\hat{\chi})| = \frac{|G|/|N|}{|\ker(\chi)|/|N|} = \frac{|G|}{|\ker(\chi)|},$ so $\mathrm{cod}(\hat{\chi}) = \mathrm{cod}(\chi)$ and therefore $\mathrm{cod}(G/N) \subseteq \mathrm{cod}(G).$
\end{proof}

\begin{Lemma}\label{subset}
    Let $G$ be a finite group with normal subgroups $N$ and $M$ such that $N\leq M$. Then, $\mathrm{cod}(G/M)\subseteq \mathrm{cod}(G/N)$.
\end{Lemma}
\begin{proof}
    By the Third Isomorphism Theorem, we know that $G / M \cong (G/N)/(M/N)$ is a quotient of  $G / N$, and by Lemma \ref{codegreeSubset}, $\mathrm{cod}(G / M) \subseteq \mathrm{cod}(G / N)$.
\end{proof}

\begin{Lemma}\label{minelements}
    Let $G$ and $H$ be finite groups such that $\mathrm{cod} (G) \subseteq \mathrm{cod} (H)$. Then there are at least $|\mathrm{cod} (G)|$ elements in $\mathrm{cod} (H)$ which divide $|G|$.
\end{Lemma}
\begin{proof}
    Let $x\in \mathrm{cod}(G)$, it is clear that $x$ divides $|G|$. Since this is true for each $x$, the lemma follows.
\end{proof}

%Further, $x=\frac{|G:\mathrm{ker}(\chi)|}{\chi(1)}$ where $\chi$ is some irreducible character of $G$.

\end{section}

\begin{section}{Main Results}

\begin{Theorem}\label{GmodN}
    Let $H$ be a sporadic simple group and $G$ a finite group with $\mathrm{cod}(G)=\mathrm{cod}(H)$. If $N$ is a maximal normal subgroup of $G$, then $G/N \cong H$.
\end{Theorem}
\begin{proof}
    By Lemma \ref{perfect}, $G$ is perfect. Thus $G/N$ is a nonabelian simple group. By Lemma \ref{subset}, we have $\mathrm{cod}(G/N) \subseteq \mathrm{cod}(G)=\mathrm{cod}(H)$. We will prove that this cannot occur unless $G/N \cong H.$ We can easily check that $\mathrm{cod}(K) \not\subseteq \mathrm{cod}(H)$ for any two non-isomorphic sporadic groups $H$ and $K$. Thus $G/N$ must belong to one of the $17$ infinite families of nonabelian simple groups. We compute the orders of these simple groups using the  formulas given in \cite{Carter}.

    For each sporadic group $H$, we check all of the possibilities for $G/N$ and restrict only to those which satisfy Lemma \ref{order}. For example, let $H\cong M$ (where $M$ denotes the Monster group). Then, we check through all $17$ families of nonabelian simple groups and return those whose orders divide $|H|$. The result is given in Table \ref{table1}. The number in the third column represents the maximal prime power $q=p^k$ which satisfies this condition (all smaller prime powers except those specified as being excluded also satisfy the condition). %The notation is standard from the ATLAS \cite{Atlas}.

    \pagebreak

    \begin{table}
    \caption{Possibilities for $G/N$ given $H \cong M$ after applying Lemma \ref{order}: $|G/N| \mid |H|.$}
    \begin{subtable}
    \centering
    \begin{tabular}{|c|c|c|}
    \hline
     $G/N$ & $n$ & Max $q=p^k$ \\
    \hline
    \hline
    $\mathrm{A}_n$ & $5$-$32$ & n/a \\
    \hline
     & $1$ & $3^4,$ excl. $2^4,2^5,43,2^6$ \\
    & $2$ & $5^2,$ excl. $11,13,2^4,19,23$ \\
    $\mathrm{L}_{n+1}(q)$ & $3$ & $3^2$ \\
    & $4$ & $2^2$ \\
    & $5$ &  $2^2$ \\
    \hline
     & $2$ &  $3^2$ \\

    & $3$  & $5$ \\

    $\mathrm{O}_{2n+1}(q)$ & $4$  & $3$ \\

     & $5$  & $2$  \\

     & $6$  & $2$  \\
    \hline
     & $3$ & $5$  \\
    $\mathrm{S}_{2n}(q)$ & $4$  & $3$  \\
     & $5$  & $2$  \\
     & $6$ & $2$  \\
    % HERE
    \hline
    \end{tabular}
    \end{subtable}
    \begin{subtable}
    \centering
    \begin{tabular}{|c|c|c|}
    \hline
     $G/N$ & $n$ & Max $q=p^k$ \\
    \hline
    \hline
     & $4$  & $3$  \\

    $\mathrm{O}_{2n}^+(q)$ & $5$  & $3$  \\

    & $6$  & $2$  \\
    \hline
    & 2 & $2^3$  \\

    $\mathrm{U}_{n+1}(q)$ & 3 & $3$  \\

    & 4 & $2$  \\

    & 5 & $2$  \\
    \hline
    $\mathrm{O}_{2n}^-(q)$ & 5 & $2$  \\

    & 6 & $2$  \\
    \hline
    $^2\mathrm{E}_{6}(q)$ & n/a & $2$  \\
    \hline
    $^3\mathrm{D}_{4}(q)$ & n/a &  $2$ \\
    \hline
    $\mathrm{Sz}(q)$ & n/a & $q=2^3,2^5$  \\
    \hline
    $^2\mathrm{F}_{4}(2)'$ & n/a & n/a  \\
    \hline
    \end{tabular}
    \end{subtable}
    \label{table1}
    \end{table}

    Now, we check each group in this table to see if it satisfies Lemma \ref{minelements}, i.e. we check how many elements in the codegree set of $H$ divide the order of $G/N$. For $H\cong M$, we find that none of the possible groups $G/N$ have order divisible by more than $3$ of the codegrees of $M$. This contradicts \cite{Alizadeh}, which shows that for any nonabelian simple group, $|\mathrm{cod}(G/N)|>3$. Thus, if $H\cong M$ and $\mathrm{cod}(G/N)\subseteq \mathrm{cod}(H),$ then $G/N\cong H$.

    We repeat this process for all of the other sporadic simple groups. For each sporadic group, $H,$ we first check which nonabelian simple groups, $G/N,$ satisfy $|G/N|$ divides $|H|.$ Second we check which of these possibilities have order divisible by more than $3$ codegrees of $H.$ We find two groups $H$ such that the number of codegrees of $H$ dividing $|G/N|$ is more than $3$. These are $H\cong\mathrm{Suz}$ with $G/N\cong O_8^+(2)$ and $H\cong\mathrm{Fi}_{23}$ with $G/N\cong O_8^+(3)$. We find $5$ and $4$ such codegrees respectively. For each of these cases, however, \cite{Aziziheris} shows that $|\mathrm{cod}(G/N)|>20$. Thus by Lemma \ref{minelements}, we cannot have $\mathrm{cod}(G/N) \subseteq \mathrm{cod}(H)$ if $G/N\not \cong H.$
\end{proof}

Now we present the proof of Theorem \ref{thm11}.%: Let $H$ be sporadic group and $G$ a finite group. If $\mathrm{cod}(G)= \mathrm{cod}(H)$, then $G \cong H$.

\begin{proof}

%We assume that the theorem is not true in search of a contradiction.
Let $G$ be a minimal counterexample and $N$ a maximal normal subgroup of $G$. By Lemma \ref{perfect}, $G$ is perfect, and by Theorem \ref{GmodN}, $G/N \cong H$. In particular, $N \neq 1$ as $G \not\cong H$.

\noindent
{\bf Step 1:} $N$ is a minimal normal subgroup of $G$.

Suppose $L$ is a nontrivial normal subgroup of $G$ with $L < N$. Then by Lemma \ref{subset}, we have $\mathrm{cod}(G/N) \subseteq \mathrm{cod}(G/L) \subseteq \mathrm{cod}(G)$. However, $\mathrm{cod}(G/N)=\mathrm{cod}(H)=\mathrm{cod}(G)$ so equality must be obtained in each inclusion. Thus, $\mathrm{cod}(G/L)=\mathrm{cod}(H)$ which implies that $G/L \cong H,$ since $G$ is a minimal counterexample. This is a contradiction since we also have $G/N\cong H,$ but $L < N$.

\noindent
{\bf Step 2:} $N$ is the only nontrivial, proper normal subgroup of $G$.

Otherwise we assume $U$ is another proper nontrivial  normal subgroup of $G$. If $N$ is included in $U$, then $U=N$ or $U=G$ since $G/N$ is simple, a contradiction. Then $N\cap U=1$ and $G=N\times U$. Since  $U$ is also a maximal normal subgroup of $G$, we have $N\cong U\cong H$. Choose $\psi_1\in \irr(N)$ and $\psi_2\in \irr(U)$ such that $\cod(\psi_1)=\cod(\psi_2)=\max(\cod(H))$. Set $\chi=\psi_1\cdot\psi_2\in \irr(G)$. Then $\cod(\chi)=(\max(\cod(H)))^2\notin \cod(G)$, a contradiction.

%To observe this, we show by contradiction that $N$ is the unique maximal normal subgroup of $G$. Suppose that $M$ is also a maximal normal subgroup of $G$. Then, $NM$ is a normal subgroup of $G$ which contains $N.$ Since $N$ is maximal, either $NM = N$ or $NM = G$. If $NM = N$, $M \leq N$ so $N=M$ since $N$ is a minimal normal subgroup of $G$. If $NM = G,$ since $N$ and $M$ are maximal normal subgroups of $G$, we can apply Theorem \ref{GmodN} to obtain $G/N\cong G/M \cong H,$ but $G/N\cong M$ and $G/M \cong N$ so $M \cong N \cong H$. Hence, $G\cong HH=H$ which contradicts $G\not \cong H.$ Thus, $N$ is the only normal subgroup of $G$ since it is minimal, maximal, and unique.

% Then if $H$ is a sporadic simple group, we can choose $\chi\in \mathrm{Irr}(H)$ such that $\mathrm{cod}(\chi)$ is maximal among $\mathrm{cod}(H)$. Then, since $G=NM \cong H, \chi^2 \in \mathrm{Irr}(G)$ and $\mathrm{cod}(\chi)^2\in \mathrm{cod}(G)$ a contradiction.

\noindent
{\bf Step 3:} For each nontrivial $\chi \in \mathrm{Irr}(G|N):=\mathrm{Irr}(G)-\mathrm{Irr}(G/N), \chi$ is faithful.

By \cite{Isaacs}*{Lemma 2.22}, we have that $\mathrm{Irr}(G/N) = \{\chi \in \mathrm{Irr}(G) |  N \leq \mathrm{ker}(\chi)\}.$ Then it follows by the definition of $\mathrm{Irr}(G|N)$ that if $\chi \in \mathrm{Irr}(G|N),$ $N \not\leq \mathrm{ker}(\chi).$ Thus since $N$ is the unique nontrivial, proper, normal subgroup of $G$, $\mathrm{ker}(\chi) = G$ or $\mathrm{ker}(\chi) = 1$. Therefore, $\mathrm{ker}(\chi) = 1$ for all nontrivial $\chi \in \mathrm{Irr}(G|N).$

\noindent
{\bf Step 4:} $N$ is an elementary abelian group.

Suppose that $N$ is not abelian. Since $N$ is a minimal normal subgroup, by \cite{Dixon}*{Theorem 4.3A (iii)}, $N=S^n$ where $S$ is a nonabelian simple group and $n\in \mathbb{Z}^+$. By Lemmas \ref{oneextend} and \ref{prodsimple}, there is a non-trivial character $\chi\in \mathrm{Irr}(N)$ which extends to some $\psi\in\mathrm{Irr}(G).$ Now, ker$(\psi)=1$ by Step 3, so cod$(\psi)=|G|/\psi(1)=|G/N|\cdot |N|/\chi(1).$ This contradicts the fact that $|G/N|$ is divisible by cod$(\psi)$, as $\chi(1)<|N|$, so $N$ must be abelian. Now to show that $N$ is elementary abelian, let a prime $p$ divide $|N|.$ Then $N$ has a $p$-Sylow subgroup $K$, and $K$ is the unique $p$-Sylow subgroup of $N$ since $N$ is abelian, so $K$ is characteristic in $N$. Thus, $K$ is a normal subgroup of $G,$ so $K=N$ as $N$ is minimal, so $|N|=p^n.$ Now, take the subgroup $N^p=\{n^p:n\in N\}$ of $N,$ which is proper by Cauchy's theorem. Then since $N^p$ is characteristic in $N,$ it must be normal in $G,$ so $N^p$ is trivial by the uniqueness of $N.$ Therefore, every element of $N$ has order $p,$ so $N$ is elementary abelian.

\noindent
{\bf Step 5:} $\mathbf{C}_G(N) = N.$

First note that since $N$ is normal, $\mathbf{C}_G(N) \trianglelefteq G.$ Additionally, since $N$ is abelian by Step 4, $N \leq \mathbf{C}_G(N),$ so by the maximality of $N,$ we must have $\mathbf{C}_G(N) = N$ or $\mathbf{C}_G(N) = G.$ If $\mathbf{C}_G(N) = N,$ we are done.

If not, then $\mathbf{C}_G(N) = G.$ Therefore $N$ must be in the center of $G.$ Then since $N$ is the unique minimal normal subgroup of $G$ by Step 2, we must have that $|N|$ is prime. If not, there always exists a proper non-trivial subgroup $K$ of $N,$ and $K$ is normal since it is contained in $\mathbf{Z}(G),$ contradicting the minimality of $N.$ Moreover, since $G$ is perfect, we have that $\mathbf{Z}(G) = N,$ and $N$ is isomorphic to a subgroup of the Schur multiplier of $G/N$ \cite{Isaacs}*{Corollary 11.20}.

If $H$ is isomorphic to any of $\mathrm{M}_{11},\mathrm{M}_{23},\mathrm{M}_{24}, \mathrm{J}_1,\mathrm{J}_4,\mathrm{Co}_2,\mathrm{Co}_3,\mathrm{Fi}_{22},\mathrm{Fi}_{23}$, He, HN, Ly, Th, or M, then by \cite{Atlas}, the Schur multiplier of $H$ is trivial, so $N=1,$ a contradiction.

If $H$ is isomorphic to $\mathrm{Co}_1,$ then $G \cong 2.\mathrm{Co}_1$ by \cite{Atlas}. But $2.\mathrm{Co}_1$ has a character degree of $24,$ which gives a codegree of $2^{19}\cdot 3^8 \cdot 5^4 \cdot 7^2 \cdot 11 \cdot 13 \cdot 23 \in \mathrm{cod}(G),$ a contradiction, since $2^{19}\cdot 3^8 \cdot 5^4 \cdot 7^2 \cdot 11 \cdot 13 \cdot 23 \notin \mathrm{cod}(H).$ If $H$ is isomorphic to $\mathrm{Fi}_{22},$ then $G \cong 2.\mathrm{Fi}_{22}$ or $G \cong 3.\mathrm{Fi}_{22}$ by \cite{Atlas}. If $G \cong 2.\mathrm{Fi}_{22},$ then $2^{13}\cdot 3^9 \cdot 5^2 \cdot 7 \cdot 13 \in \mathrm{cod}(G),$ a contradiction. If $G \cong 3.\mathrm{Fi}_{22},$ then $2^{17}\cdot 3^7 \cdot 5^2 \cdot 6 \cdot 11 \in \mathrm{cod}(G),$ a contradiction.

Similarly, for any sporadic simple group $H$ with non-trivial Schur multiplier, we use \cite{Atlas} to reach a contradiction as we did above, by finding an element of $\mathrm{cod}(G)$ that is not in $\mathrm{cod}(H).$

Thus $\mathbf{C}_G(N) = N.$

\noindent
{\bf Step 6:} Let $\lambda$ be a non-trivial character in $\mathrm{Irr}(N)$ and $\vartheta \in \mathrm{Irr}(I_G(\lambda)|\lambda),$ the set of irreducible constituents of $\lambda^{I_G(\lambda)},$ where $I_G(\lambda)$ is the inertia group of $\lambda \in G.$ Then $\frac{|I_G(\lambda)|}{\vartheta(1)} \in \mathrm{cod}(G).$ Also, $\vartheta(1)$ divides $|I_G(\lambda)/N|,$ and $|N|$ divides $|G/N|.$ Lastly, $I_G(\lambda) < G,$ i.e. $\lambda$ is not $G$-invariant.

Let $\lambda$ be a non-trivial character in $\operatorname{Irr}(N)$ and $\vartheta \in \operatorname{Irr}(I_G(\lambda)|\lambda)$. Let $\chi$ be an irreducible constituent of $\vartheta^G.$ By \cite{Isaacs}*{Corollary 5.4}, we know $\chi \in \operatorname{Irr}(G)$, and by \cite{Isaacs}*{Definition 5.1}, we have $\chi(1) = \frac{|G|}{|I_G(\lambda)|} \cdot \vartheta(1)$. Moreover, we know tat $\operatorname{ker}(\chi) = 1$ by Step 2, and thus $\cod(\chi) = \frac{|G|}{\chi(1)} = \frac{|I_G(\lambda)|}{\vartheta(1)}$, so $\frac{|I_G(\lambda)|}{\vartheta(1)} \in \mathrm{cod}(G)$. Now, since $N$ is abelian, $\lambda(1) = 1$, so we have $\vartheta(1) = \vartheta(1)/\lambda(1)$ which divides $\frac{|I_G(\lambda)|}{|N|}$, so $|N|$ divides $\frac{|I_G(\lambda)|}{\vartheta(1)}$. Moreover, we know that $\cod(G) = \cod(G/N),$ and all elements in $\cod(G/N)$ divide $|G/N|$, so $|N|$ divides $|G/N|$.

Next, we want to show $I_G(\lambda)$ is a proper subgroup of $G$. To reach a contradiction, assume $I_G(\lambda) = G$. Then $\operatorname{ker}(\lambda) \unlhd G$. From Step 2, we know $\operatorname{ker}(\lambda) = 1,$ and from Step 4, we know $N$ is a cyclic group of prime order. Thus by the Normalizer-Centralizer theorem, we have $G / N= \mathbf{N}_{G}(N) / \mathbf{C}_{G}(N) \leq \operatorname{Aut}(N)$ so $G / N$ is abelian, a contradiction.
% By Theorem \ref{clifford}, since $\chi \cong \theta$ we'll have $e = 1$, and $\theta^g(1) = \theta(1)$ by definition. Thus we conclude $\chi(1) = |G : I_G(\lambda)| \cdot \theta(1)$.

\noindent
{\bf Step 7:} Final contradiction.

From Step 4, $N$ is an elementary abelian group of order $p^n$ for some prime $p$ and integer $n\geq1$. By the Normalizer-Centralizer theorem, $H \cong G/N = \mathbf{N}_G(N)/\mathbf{C}_G(N) \leq \mathrm{Aut}(N)$ and $n>1$. Note that in general, $\mathrm{Aut}(N)=\mathrm{GL}(n,p)$. By Step 6, $|N|$ divides $|G/N|,$ so we only need to consider primes $p$ such that $p^2$ divides $|H|$. For instance, if $H\cong M_{11}$, $|H|=2^4\cdot3^2\cdot5\cdot11$ so $|N|=2^2, 2^3, 2^4$  or $3^2$. Then, we can compute $|\mathrm{Aut}(N)|=|\mathrm{GL}(n,p)|$ which is equal to $6, 168, 20160,$ and $48$ respectively. In each of these four cases, $|H|\nmid  |\mathrm{Aut}(N)|$.

For each sporadic group $H$, we follow a similar procedure to computationally check which possibilities of $(p,n)$ satisfy $p^n$ divides $|H|$ and $|H|$ divides $|\mathrm{GL}(n,p)|$. We find only the following $7$ possible $H$ which satisfy this condition, listed in Table \ref{tblexc}. In other words, Table \ref{tblexc} gives all groups $H$ and pairs $(p,n)$ such that $p^n$ divides $|H|$ and $|H|$ divides $|\mathrm{GL}(n,p)|$.

\begin{table}[H]
    \centering
    \caption{Minimum degree of faithful representations of sporadic groups over $\mathbb{F}_{p^k}.$}
    \begin{tabular}{ |c|c|c|c| }
    \hline
    Group & $p$ & $n$ & Minimum Degree \\
    \hline
    \hline
    He & $2$ & $9-10$ & $51$ \\
    \hline
    Suz & $2$ & $12-13$ & $110$ \\
    \hline
    $\mathrm{Fi}_{22}$ & $2,3$ & $14-17,$ $8-9$ & $78,77$ \\
    \hline
    $\mathrm{Fi}_{23}$ & $2$ & $18$ & $782$ \\
     \hline
    $\mathrm{Co}_{2}$ & $2$ & $12-18$ & $22$ \\ \hline
    $\mathrm{Co}_{1}$ & $2$ & $16-21$ & $24$ \\ \hline
    B & $2$ & $23-41$ & $4370$ \\
    \hline
    \end{tabular}
    \label{tblexc}
\end{table}

The final column in the table above gives the minimal degree of a faithful representation of the group $H$ over a finite field of characteristic $p$  \cite{Jansen}. As this minimal degree is larger than the largest possible $n$ in each case, we deduce that $H$ cannot be isomorphic to a subgroup of $\mathrm{GL}(n,p)$. Hence, we have a contradiction for any sporadic simple group $H$. Thus, $N=1$ and $G \cong H$.
\end{proof}

\end{section}

\begin{section}{Acknowledgements}

This research was conducted under NSF-REU grant DMS-1757233, DMS-2150205 and NSA grant H98230-21-1-0333, H98230-22-1-0022 by Dolorfino, Martin, Slonim, and Sun during the Summer of 2022 under the supervision of Yang. The authors gratefully acknowledge the financial support of NSF and NSA, and also thank Texas State University for providing a great working environment and support. Yang was also partially supported by grants from the Simons Foundation (\#499532, \#918096, YY).\\
\end{section}

\noindent \textbf{Competing interests} The authors declare none.

\bigskip
\noindent
\textbf{Data availability Statement:} Data sharing not applicable to this article as no datasets were generated or analysed during the current study.

\vspace{-10px}
\end{document}